\documentclass[12pt,leqno]{article}
\usepackage{soul}
\usepackage[doc]{optional}
\usepackage{xcolor}
\definecolor{labelkey}{rgb}{0,0.08,0.45}
\definecolor{rekey}{rgb}{0,0.6,0.0}
\definecolor{Brown}{rgb}{0.45,0.0,0.05}
\usepackage{exscale,relsize}
\usepackage{amsmath}
\usepackage{amsfonts}
\usepackage{amssymb}
\usepackage{calc}
\usepackage{theorem}
\usepackage{pifont}      
\usepackage{graphicx}
\newcommand{\scal}[2]{\langle{{#1},{#2}}\rangle}

\newcommand{\RR}{\ensuremath{\mathbb R}}

\newcommand{\NN}{\ensuremath{\mathbb N}}

\newcommand{\menge}[2]{\big\{{#1} \mid {#2}\big\}}

\newcommand{\To}{\ensuremath{\rightrightarrows}}

\newcommand{\spand}{\operatorname{span}}

\newcommand{\dom}{\ensuremath{\operatorname{dom}}}

\newcommand{\gra}{\ensuremath{\operatorname{gra}}}

\newcommand{\ran}{\ensuremath{\operatorname{ran}}}

\newcommand{\rank}{\ensuremath{\operatorname{rank}}}

\newcommand{\Id}{\ensuremath{\operatorname{Id}}}

\renewcommand{\phi}{\ensuremath{\varphi}}

\newcommand{\bA}{\ensuremath{\mathbf{A}}}

\newcommand{\vp}{\ensuremath{\mathbf{V_{+}}}}
\newcommand{\vm}{\ensuremath{\mathbf{V_{-}}}}
\newcommand{\vn}{\ensuremath{\mathbf{V_{0}}}}
\newcommand{\wtg}{\ensuremath{\widetilde{G}}}

\newcommand{\whg}{\ensuremath{\widehat{G}}}

\newtheorem{theorem}{Theorem}[section]
\newtheorem{lemma}[theorem]{Lemma}
\newtheorem{fact}[theorem]{Fact}
\newtheorem{corollary}[theorem]{Corollary}
\newtheorem{proposition}[theorem]{Proposition}

\theoremstyle{plain}{\theorembodyfont{\rmfamily}
}
\theoremstyle{plain}{\theorembodyfont{\rmfamily}
}
\theoremstyle{plain}{\theorembodyfont{\rmfamily}
}
\theoremstyle{plain}{\theorembodyfont{\rmfamily}
\newtheorem{example}[theorem]{Example}}
\theoremstyle{plain}{\theorembodyfont{\rmfamily}
\newtheorem{remark}[theorem]{Remark}}
\theoremstyle{plain}{\theorembodyfont{\rmfamily}
}

\begin{document}

\title{{\sffamily  Maximally Monotone Linear Subspace Extensions of Monotone Subspaces:
Explicit Constructions and Characterizations}}

\author{Xianfu Wang
\thanks{Mathematics, Irving K.\ Barber School,
The University of British Columbia Okanagan, Kelowna,
B.C. V1V 1V7, Canada.
Email:
\texttt{shawn.wang@ubc.ca}.} ~and~ Liangjin Yao
\thanks{Mathematics,
Irving K.\ Barber School,
The University of British Columbia Okanagan, Kelowna,
B.C. V1V 1V7, Canada.
E-mail:  \texttt{ljinyao@interchange.ubc.ca}.}
\\
Dedicated to Jonathan Borwein on the occasion of his 60th birthday.}

 \vskip 3mm

\date{March 6, 2011}
\maketitle
\begin{abstract} \noindent
Monotone linear relations play important roles in variational inequality problems and quadratic optimizations.
In this paper, we give explicit maximally monotone linear subspace extensions of a monotone linear relation in finite
dimensional spaces. Examples are provided to illustrate our extensions. Our results generalize a recent result by
Crouzeix and Anaya.
\end{abstract}



\section{Introduction}
Throughout this paper, we assume that
$\RR^n$ ( $n\in\NN=\{1,2,3,\ldots\}$) is an Euclidean space with the
inner product $\scal{\cdot}{\cdot}$,
and induced Euclidean norm $\|\cdot\|$.
Let $G\colon \RR^n\To\RR^n$
be a \emph{set-valued operator}
from $\RR^n$ to $\RR^n$, i.e., for every $x\in \RR^n$, $Gx\subseteq \RR^n$,
and let
$\gra G = \menge{(x,x^*)\in \RR^n\times \RR^n}{x^*\in Gx}$ be
the \emph{graph} of $G$.
Recall that $G$ is  \emph{monotone} if
\begin{equation}
\big(\forall (x,x^*)\in \gra G\big)\big(\forall (y,y^*)\in\gra
G\big) \quad \scal{x-y}{x^*-y^*}\geq 0,
\end{equation}
and \emph{maximally monotone} if $G$ is monotone and $G$ has no proper monotone extension
(in the sense of graph inclusion).
We say that $G$ is a \emph{linear relation} if $\gra G$ is a linear subspace.
While linear relations have been extensively studies \cite{Cross,Yao1,borlewis,bwy09,BWY8,voisei}, monotone operators
are ubiquitous in convex optimization and variational analysis \cite{BC2011,RockWets,borlewis,BorVan}.

The central object of this paper is to consider the linear relation
$G:\RR^n\To\RR^n$:
\begin{enumerate}
\item $\gra G=\{(x,x^*)\mid Ax+Bx^*=0\}.$
\item  $A,B\in\RR^{p\times n}$.
\item $\rank(A\ B)=p$.
\end{enumerate}

Our main concern is to find explicit maximally monotone linear subspace extensions of $G$.
Recently, finding constructive maximal monotone extensions instead of
using Zorn's lemma has been a very active
topic \cite{BW,bw10,CrouAna1,CrouAna2,CrouAna}.
In \cite{CrouAna}, Crouzeix and Anaya gave an algorithm to find maximally monotone linear subspace
extensions of $G$, but it is not clear what the maximally monotone extensions are analytically. In this paper, we provide some maximally monotone extensions
of $G$ with closed analytical forms. Along the way, we also give a new proof to Crouzeix and Anaya's characterizations on monotonicity and maximal monotonicity of $G$. Our key tool is the Brezis-Browder characterization
of maximally monotone linear relations.

The paper is organized as follows.
In the remainder of this introductory section, we describe some
central notions fundamental to our analysis.
In Section~\ref{s:aux}, we collect some auxiliary results for future reference
and for the
reader's convenience.
Section~\ref{s:main} provides explicit self-dual maximal monotone extensions by using subspaces on which
$AB^{\intercal}+BA^{\intercal}$ is negative semidefinite, and obtain a complete characterization of
all maximal monotone extensions.
Section~\ref{Minty} deals with  Minty's parameterizations of monotone operator $G$. In Section~\ref{DRan:1},
we get some explicit maximally monotone extensions with the same domain or the same range by utilizing normal cone
operators.
In Section~\ref{Exam:sec}, we illustrate our maximally monotone extensions by considering three examples.

Our notations are standard.
We use $\dom G = \menge{x\in \RR^n}{Gx\neq\varnothing}$ for the \emph{domain} of $G$,
$\ran G=G(\RR^n)$ for the \emph{range} of $G$ and
$\ker G=\menge{x\in \RR^n}{0\in Gx}$ for the \emph{kernal} of $G$.
Given a subset $C$ of $\RR^n$,
$\spand C$ is the \emph{span} (the set of all finite linear
combinations) of $C$. We set $$C^{\bot}=
\{x^*\in \RR^n \mid(\forall c\in C)\, \langle x^*, c\rangle=0\}.$$
Then the \emph{adjoint} of $G$, denoted by $G^*$, is defined
by
\begin{equation*}
\gra G^* =
\{(x,x^*)\in \RR^n\times \RR^n\mid (x^*,-x)\in (\gra G)^\bot\}.
\end{equation*}
The set $\RR^{n\times p}$ is the set of all the $n\times p$ matrices, for $n,p\in \NN$.
Then $\rank(M)$ is the \emph{rank} of the matrix $M\in\RR^{n\times p}$.
Let $\Id:\RR^n\rightarrow \RR^n$ denote the identity mapping, i.e., $\Id x=x$ for $x\in \RR^n$.
We also set $P_X: \RR^n\times \RR^n\rightarrow \RR^n\colon (x,x^*)\mapsto x$,
and $P_{X^*}: \RR^n\times \RR^n\rightarrow \RR^n\colon (x,x^*)\mapsto x^*$.
If $\mathcal{X}, \mathcal{Y}$ are subspaces of $\RR^n$, we let
$$\mathcal{X}+\mathcal{Y}=\menge{x+y}{x\in \mathcal{X}, y\in\mathcal{Y}}.$$

 Counting multiplicities, let
\begin{equation}\label{t:pos}
\lambda_{1},\lambda_2,\cdots,\lambda_k \text{  be
all positive eigenvalues of $(AB^{\intercal}+BA^{\intercal})$ and}
\end{equation}
\begin{equation}\label{t:neg}
 \lambda_{k+1},\lambda_{k+2},\cdots,\lambda_p \text{ be nonpositive eigenvalues of $(AB^{\intercal}+BA^{\intercal})$.}
 \end{equation}
Moreover, let $v_i$ be an  eigenvector of eigenvalue $\lambda_i$ of $(AB^{\intercal}+BA^{\intercal})$ satisfying
$\|v_i\|=1$, and $\langle v_i,v_j\rangle=0$ for $ \ 1\leq i\neq j\leq q$.
It will be convenient to put
\begin{align}\label{t:vmatrix}
\Id_{\lambda}=\text{diag}(\lambda_{1},\cdots, \lambda_{p})=
\begin{pmatrix}
\lambda_1 & 0 & 0 & \cdots  &  0 \\
0 & \lambda_2 & 0 & \cdots &   0 \\
0 & 0 & \lambda_3 &           & \vdots \\
\vdots & 0 & 0      & \ddots &   0 \\
0 & 0 & 0 & 0 & \lambda_p \\
\end{pmatrix},
\quad V=\left[v_1\,v_2\,\cdots\,v_p\right].
\end{align}

\section{Auxiliary results on linear relations}\label{s:aux}

In this section, we collect some facts and preliminary results which will be used in sequel.

We first provide a result about subspaces on which a linear operator from $\RR^n\to \RR^n$, i.e,
an $n\times n$ matrix, is monotone.
For $M\in \RR^{n\times n}$, define three subspaces of $\RR^n$, namely, the
 positive eigenspace, null eigenspace and negative eigenspace associated with $M+M^{\intercal}$
  by
$$
\vp(M) =\spand\begin{cases} w_{1},\cdots, w_{s}: &\ w_{i} \text{ is an eigenvector of positive eigenvalue $\alpha_{i}$ of $M+M^{\intercal}$}\\
& \langle w_{i}, w_{j}\rangle =0 \ \forall \ i\neq j, \|w_{i}\|=1,  i,j=1,\ldots, s.
\end{cases}
\Bigg \}
$$
$$
\vn(M) =\spand\begin{cases} w_{s+1},\cdots, w_{l}: &\ w_{i} \text{ is an eigenvector of $0$ eigenvalue of $M+M^{\intercal}$}\\
& \langle w_{i}, w_{j}\rangle =0, \ \forall \ i\neq j, \|w_{i}\|=1 , i,j=s+1,\ldots, l.
\end{cases}
\Bigg \}
$$
$$
\vm(M) =\spand\begin{cases} w_{l+1},\cdots, w_{n}: &\ w_{i} \text{ is an eigenvector of negative eigenvalue $\alpha_{i}$ of $M+M^{\intercal}$}\\
& \langle w_{i}, w_{j}\rangle =0 \ \forall \ i\neq j, \|w_{i}\|=1, i,j=l+1,\ldots, n.
\end{cases}
\Bigg \}
$$
which is possible since a symmetric matrix always has a complete orthonormal set of eigenvectors,
\cite[pages 547--549]{Meyer}.
\begin{proposition}\label{t:decomposition}
Let $M$ be an $n\times n$ matrix. Then
\begin{enumerate}
\item $M$ is strictly monotone on $\vp(M)$. Moreover,
$M+M^{\intercal}:\vp(M)\rightarrow \vp(M)$ is a bijection.

\item $M$ is monotone on $\vp(M)+ \vn(M)$.
\item $-M$ is strictly monotone on $\vm(M)$. Moreover,
$-(M+M^{\intercal}):\vm(M)\rightarrow \vm(M)$ is a bijection.

\item $-M$ is monotone on $\vm(M)+\vn(M)$.
\item For every $x\in \vn(M)$, $(M+M^{\intercal})x=0$ and $\langle x,Mx\rangle=0$.
\end{enumerate}
In particular, the orthogonal decomposition holds: $\RR^n=\vp(M)\oplus\vn(M)\oplus\vm(M)$.
\end{proposition}
\begin{proof}
 (i):  Let $x\in \vp(M)$. Then $x=\sum_{i=1}^{s}l_{i}w_{i}$ for some $(l_{1},\ldots, l_{s})\in \RR^s$. Since
$\{w_{1},\cdots, w_{s}\}$ is a set of orthonormal vectors, they are linearly independent
so that
$$x\neq 0\quad \Leftrightarrow\quad (l_{1},\cdots, l_{s})\neq 0.$$
Note that $\alpha_{i}>0$ when $i=1,\ldots, s$ and $\langle w_{i}, w_{j}\rangle=0$ for $i\neq j$.
We have
\begin{align*}
2\langle x, Mx\rangle & = \langle x, (M+M^{\intercal}) x\rangle =\langle\sum_{i=1}^{s}
l_{i}w_{i}, (M+M^{\intercal})( \sum_{i=1}^{s}
l_{i}w_{i})\rangle\\
& = \langle \sum_{i=1}^{s}
l_{i}w_{i}, \sum_{i=1}^{s}
l_{i}\alpha_{i}w_{i}\rangle =\sum_{i=1}^{s}\alpha_{i}l_{i}^2>0
\end{align*}
if $x\neq 0$.

For every $x\in \vp(M)$ with $x=\sum_{i=1}^{s}l_{i}w_{i}$, we have
$$(M+M^{\intercal})x=\sum_{i=1}^{s}l_{i}(M+M^{\intercal})w_{i}=\sum_{i=1}^{s}\alpha_{i}l_{i}w_{i}\in\vp(M).$$
As $\alpha_{i}>0$ for $i=1,\ldots,s$ and $\{w_{1},\ldots, w_{s}\}$ is an orthonormal basis of $\vp(M)$, we conclude that
$M+M^{\intercal}: \vp(M)\rightarrow \vp(M)$ is a bijection.

(ii): Let $x\in \vp(M)+\vn(M)$. Then $x=\sum_{i=1}^{l}l_{i}w_{i}$ for some $(l_{1},\ldots, l_{l})\in \RR^l$.
Note that $\alpha_{i}\geq 0$ when $i=1,\ldots, l$ and $\langle w_{i}, w_{j}\rangle=0$ for $i\neq j$.
We have
\begin{align*}
2\langle x, Mx\rangle & = \langle x, (M+M^{\intercal}) x\rangle =\langle\sum_{i=1}^{l}
l_{i}w_{i}, (M+M^{\intercal})( \sum_{i=1}^{l}
l_{i}w_{i})\rangle\\
& = \langle \sum_{i=1}^{l}
l_{i}w_{i}, \sum_{i=1}^{l}
l_{i}\alpha_{i}w_{i}\rangle =\sum_{i=1}^{l}\alpha_{i}l_{i}^2\geq 0.
\end{align*}
The proofs for (iii), (iv) are similar as (i), (ii).

(v): For $x\in \vn(M)$,
$$2\langle x,Mx\rangle=\langle x, (M+M^{\intercal})x\rangle=\langle x, 0\rangle=0.$$
\end{proof}
\begin{corollary}\label{t:split}
 Then following holds.
\begin{enumerate}
\item
$$\gra T= \{(B^{\intercal}u, A^{\intercal}u)\mid u\in \vp(BA^{\intercal})\}$$
is strictly monotone.
\item
$$\gra T= \{(B^{\intercal}u, A^{\intercal}u)\mid u\in \vp(BA^{\intercal})+\vn(BA^{\intercal})\}$$
is monotone.
\item $$\gra T= \{(B^{\intercal}u, -A^{\intercal}u)\mid u\in \vm(BA^{\intercal})\}$$
is strictly monotone.
\item $$\gra T= \{(B^{\intercal}u, -A^{\intercal}u)\mid u\in \vm(BA^{\intercal})+\vn(BA^{\intercal}))\}$$
is monotone.
\end{enumerate}
\end{corollary}
\begin{proof}
It follows from Proposition~\ref{t:decomposition} and $\langle B^{\intercal}u, A^{\intercal}u\rangle=\langle u, BA^{\intercal}u\rangle,\ \forall u\in\RR^n.$
\end{proof}
\begin{lemma} \label{t:dimensionnochange}
For every subspace $S\subseteq \RR^p$, the following hold.
\begin{equation}\label{t:dim+}
\dim\{(B^{\intercal}u, A^{\intercal}u)\mid u\in S\}=\dim S.
\end{equation}
\begin{equation}\label{t:dim-}
\dim\{(B^{\intercal}u, -A^{\intercal}u)\mid u\in S\}=\dim S.
\end{equation}
\end{lemma}
\begin{proof} Let $\dim S=t$ and $\{u_{1},\ldots, u_{t}\}$ be a basis of $S$.
We claim that the set of vectors
$$\left\{\left(\begin{matrix} B^{\intercal} \\
A^{\intercal}
\end{matrix}\right) u_{i}\mid \ i=1,\ldots,t\right\}
$$
is linearly independent. Indeed, because $(A\ B)$ has full row rank $p$,
$$\sum_{i=1}^{t}l_{i}\left(\begin{matrix} B^{\intercal} \\
A^{\intercal}
\end{matrix}\right)u_{i}=\left(\begin{matrix} B^{\intercal} \\
A^{\intercal}
\end{matrix}\right) \sum_{i=1}^{t}l_{i} u_{i}=0\quad \Leftrightarrow\quad \sum_{i=1}^{t}l_{i}u_{i}=0\quad \Leftrightarrow
\quad l_{i}=0 \text{ for $i=1,\ldots,t$}.$$
Note that
$$\left(\begin{matrix} B^{\intercal}u \\
-A^{\intercal}u
\end{matrix}\right)=\left(\begin{matrix} \Id & 0\\
0 & -\Id\end{matrix}\right)\left(\begin{matrix} B^{\intercal}u \\
A^{\intercal}u
\end{matrix}\right)
$$
and
$$\left(\begin{matrix} \Id & 0\\
0 & -\Id\end{matrix}\right)$$
is invertible,
we have
$$\dim\{(B^{\intercal}u, A^{\intercal}u)\mid u\in S\}=\dim\{(B^{\intercal}u, -A^{\intercal}u)\mid u\in S\}$$
so \eqref{t:dim-} follows from \eqref{t:dim+}. Alternatively, see \cite[page 208, Exercise 4.49]{Meyer}.
\end{proof}

\begin{fact}\label{FaAE:1} We have
\begin{align*}(AB^{\intercal}+BA^{\intercal})V=V\Id_{\lambda}.\end{align*}
\end{fact}
\begin{proof}
Let $y=\left[y_1\,y_2\,\cdots y_p\right]^{\intercal}\in\RR^p$. Then we have
\begin{align*}
(AB^{\intercal}+BA^{\intercal})Vy
=  \big(AB^{\intercal}+BA^{\intercal}\big)(\sum_{i=1}^{p}y_iv_i)
=\sum_{i=1}^{p}\lambda_i y_iv_i
=V\Id_{\lambda}y.
\end{align*}
\end{proof}

Two key criteria concerning maximally monotone linear relations come as follows:
\begin{fact}\emph{(See \cite[Proposition~4.2.9 ]{Yao} or \cite[Proposition~ 2.10]{BWY8}.)}\label{TheSIAM}
Let $T\colon \RR^n\To \RR^n$ be a monotone linear relation.
The following are equivalent:
\begin{enumerate}
\item $T$ is maximally monotone.
\item $\dim\gra T=n$.
\item $\dom T=(T0)^{\perp}.$
\end{enumerate}
\end{fact}

\begin{fact}[Br\'ezis-Browder] \emph{(See \cite[Theorem~2]{Brezis-Browder}, or \cite{Yao1} or \cite{Si3}.)}
\label{Sv:7}
Let  $T\colon \RR^n \rightrightarrows \RR^n$ be a monotone linear relation.
 Then the following statements are equivalent.
\begin{enumerate}
 \item $T$ is maximally monotone.
\item  $T^*$ is maximal  monotone.
\item   $T^*$ is monotone.
\end{enumerate}
\end{fact}

Some basic properties of $G$ are:
\begin{lemma}\label{AuSe:1}
\begin{enumerate}
\item $\gra G=\ker(A\ B)$.
\item $G0=\ker B, G^{-1}(0)=\ker A$.
\item $\dom G= P_{X}(\ker(A\ B))$ and $\ran G= P_{X^*}(\ker (A\ B))$.
\item $\ran (G+\Id)=P_{X^*}(\ker(A-B\ B))=P_{X}(\ker(A \ B-A))$, and
$$\dom G=P_{X}(\ker( A-B \ B)), \quad \ran G=P_{X^*}(\ker (A \ (B-A)).$$
\item $\dim G=2n-p$.
\end{enumerate}
\end{lemma}
\begin{proof}
(i), (ii), (iii) follow from definition of $G$.
Since
$$Ax+Bx^*=0\quad \Leftrightarrow\quad (A-B)x+B(x+x^*)=0 \quad \Leftrightarrow\quad
A(x+x^*)+(B-A)x^*=0,$$
(iv) holds.

(v): We have
\begin{align*}
2n=\dim\ker (A\ B)+\dim\ran \left(\begin{matrix}A^{\intercal}\\
B^{\intercal}\end{matrix}\right)=\dim G+p.\end{align*} Hence $\dim G=2n-p$.
\end{proof}

The following result summarizes the monotonicities of $G^*$ and $G$.
\begin{lemma}\label{MonL:1} The following hold.
\begin{enumerate}\item  $\gra G^*=\{(B^{\intercal}u,-A^{\intercal}u)\mid u\in\RR^p\}$.
\item  $G^*$ is monotone $\Leftrightarrow$ $A^{\intercal}B+B^{\intercal}A$ is negative-semidefinite.
 \item Assume $G$ is monotone. Then $n\leq p$. Moreover, $G$ is maximally monotone
 if and only $\dim G=n=p$.
\end{enumerate}
\end{lemma}
\begin{proof}
(i):\, By Lemma~\ref{AuSe:1}(i), we have
\begin{align*}
(x,x^*)\in \gra G^*\Leftrightarrow (x^*,-x)\in \gra G^{\bot}=\ran \left(\begin{matrix}A^{\intercal}\\
B^{\intercal}\end{matrix}\right)=\{(A^{\intercal}u,B^{\intercal}u)\mid u\in\RR^p\}.
\end{align*}
Thus
 $\gra G^*=\{(B^{\intercal}u,-A^{\intercal}u)\mid u\in\RR^p\}$.

 (ii):\, Since $\gra G^*$ is a linear subspace, by (i),
 \begin{align*}
& G^*\, \text{is monotone}\Leftrightarrow \langle B^{\intercal}u,-A^{\intercal}u\rangle\geq0,\quad\forall u\in\RR^p
 \Leftrightarrow \langle u,-BA^{\intercal}u\rangle\geq0,\quad\forall u\in\RR^p\\
 &\Leftrightarrow\langle u,BA^{\intercal}u\rangle\leq0,\quad\forall u\in\RR^p
 \Leftrightarrow \langle u,(A^{\intercal}B+B^{\intercal}A)u\rangle\leq0,\quad\forall u\in\RR^p\\
& \Leftrightarrow (A^{\intercal}B+B^{\intercal}A)\, \text{is negative semidefinite}.
\end{align*}

(iii):\, By Fact~\ref{TheSIAM} and Lemma~\ref{AuSe:1}(v), $2n-p=\dim \gra G\leq n\Rightarrow n\leq p$.
 By Fact~\ref{TheSIAM} and Lemma~\ref{AuSe:1}(v) again,
$G$ is maximally monotone $\Leftrightarrow$  $2n-p=\dim \gra G= n\Leftrightarrow\dim \gra G=p =n$.
\end{proof}

\section{Explicit maximal monotone extensions of monotone linear relations}\label{s:main}

In this section, we give explicit maximal monotone linear subspace extensions of $G$ by using $\vp(AB^{\intercal})$ or
$V_{g}$. A characterization of all the maximally monotone extensions of $G$ is also given.

\begin{lemma}\label{FaAEL:2} Let $N \in\RR^{p\times p}$ and  $\widetilde{G}$ and $\whg$ be defined  by
\begin{align*}\gra\widetilde{G}&=\{(x,x^*)\mid N^{\intercal}V^{\intercal}Ax+N^{\intercal}V^{\intercal}Bx^*=0\}\\
\gra \whg&=\{(B^{\intercal}u, -A^{\intercal}u)\mid\ u\in\ran VN)\}.\end{align*}
Then $(\widetilde{G})^*=\whg$.
\end{lemma}
\begin{proof}Let $(y,y^*)\in \RR^n\times \RR^n$. Then we have
\begin{align*}
(y,y^*)\in\gra (\widetilde{G})^*&\Leftrightarrow (y^*,-y)\in(\gra\widetilde{G})^{\bot}
=(\ker\left(\begin{matrix}N^{\intercal}V^{\intercal}A\quad
N^{\intercal}V^{\intercal}B\end{matrix}\right)\big)^{\bot}=\ran\left(\begin{matrix}A^{\intercal}VN\\
B^{\intercal}VN\end{matrix}\right)\\
&\Leftrightarrow (y,y^*)\in \gra \whg.
\end{align*}
Hence $(\widetilde{G})^*=\whg$.
\end{proof}

\begin{lemma}\label{wtgtog}
Define $\widetilde{G}$ and $\whg$ by
\begin{align*}\gra\widetilde{G}&=\{(x,x^*)\mid V_g Ax+V_g Bx^*=0\}\\
\gra \whg&=\{(B^{\intercal}u, -A^{\intercal}u)\mid\ u\in\vm(BA^{\intercal})+\vn(BA^{\intercal})\},\end{align*}
where $V_g$ is $(p-k)\times p$ matrix defined by
\begin{align*}V_g=\left(\begin{matrix}
v_{k+1}^{\intercal}\\
v_{k+2}^{\intercal}\\
\vdots\\
v_{p}^{\intercal}
\end{matrix}
\right).
\end{align*}

Then
\begin{enumerate}

\item $\whg$ is monotone.
\item
$
(\whg)^* =\wtg.
$

\item $\gra\widetilde{G}=\gra G+\left\{\left(\begin{matrix}B^{\intercal}\\
A^{\intercal}\end{matrix}\right)u\mid u\in \vp(BA^{\intercal})\right\}.$

\end{enumerate}
\end{lemma}

\begin{proof}
(i): Apply Corollary~\ref{t:split}(iv).

(ii): Notations are as in \eqref{t:vmatrix}.
Let
\begin{align}N=\left[0\ 0\cdots 0\ e_{k+1}\cdots e_p\right],\label{NewMatr:1}\end{align}
where $e_i=\left[0,0\cdots 1,0\cdots0\right]^{\intercal}$: the ith entry is 1 and the others are $0$.

Then we have
\begin{align}
N^{\intercal}V^{\intercal}=\bigg((v_{1}\ \cdots v_{k}\ \ V_{g}^{\intercal})
\begin{pmatrix}
\mathbf{0} & \mathbf{0}\\
\mathbf{0} & \Id
\end{pmatrix}
\bigg)^{\intercal}=
\left(\begin{matrix}0\\ V_g\end{matrix}\right).
\end{align}
Then we have
\begin{align*}
V_g Ax +V_g Bx^*&=0\Leftrightarrow \left(\begin{matrix}0\\V_g Ax +V_g Bx^*\end{matrix}\right)=0\\
&\Leftrightarrow
N^{\intercal}V^{\intercal} Ax+N^{\intercal}V^{\intercal}Bx^*=0,\quad \forall (x,x^*)\in \RR^n\times \RR^n.
\end{align*}
Hence \begin{align*}\gra\widetilde{G}&=\{(x,x^*)\mid N^{\intercal}V^{\intercal} Ax+N^{\intercal}V^{\intercal}Bx^*=0\}.
\end{align*}
Thus by Lemma~\ref{FaAEL:2},
\begin{align*}
\gra(\widetilde{G})^*=\{(B^{\intercal}u, -A^{\intercal}u)\mid u\in\ran VN=\ran\left(\begin{matrix}0\ V^{\intercal}_g\end{matrix}\right)=\vm(BA^{\intercal})+\vn(BA^{\intercal})\}=\gra \whg.
\end{align*}
Hence $(\whg)^*=(\wtg)^{**}=\wtg$.

(iii):
Let $J$ be defined by
\begin{align*}
\gra J=
\gra G+ \left\{\left(\begin{matrix}B^{\intercal}\\
A^{\intercal}\end{matrix}\right)u\mid u\in \vp(BA^{\intercal})\right\}.
\end{align*}

Then we have
$$(\gra J)^{\bot}= (\gra G)^{\perp}\cap \left\{\left(\begin{matrix}B^{\intercal}\\
A^{\intercal}\end{matrix}\right)u\mid u\in \vp(BA^{\intercal})\right\}^{\perp}.$$
By Lemma~\ref{AuSe:1}(i),
$$\gra G^{\perp}=\left\{\left(\begin{matrix}A^{\intercal}\\
B^{\intercal}\end{matrix}\right)w\mid w\in \RR^{p}\right\}
$$
Then
$$\left(\begin{matrix}A^{\intercal}\\
B^{\intercal}\end{matrix}\right)w\in \left\{\left(\begin{matrix}B^{\intercal}\\
A^{\intercal}\end{matrix}\right)u\mid u\in \vp(BA^{\intercal})\right\}^{\perp}$$
if and only if
$$\langle (A^{\intercal}w, B^{\intercal}w), (B^{\intercal}u, A^{\intercal}u)\rangle =0\quad \forall \ u\in \vp(BA^{\intercal}),$$
that is,
\begin{equation}\label{t:party}
\langle A^{\intercal} w, B^{\intercal}u\rangle +\langle B^{\intercal} w, A^{\intercal} u\rangle=
\langle w, (AB^{\intercal}+BA^{\intercal})u\rangle=0\quad
\forall u\in \vp(AB^{\intercal}).
\end{equation}
Because
$AB^{\intercal}+BA^{\intercal}: \vp(AB^{\intercal})\mapsto \vp(AB^{\intercal})$ is
 onto by Proposition~\ref{t:decomposition}(i), we obtain that \eqref{t:party}
holds if and only if
$w\in \vm(AB^{\intercal}) +\vn(AB^{\intercal})$.
Hence
$$(\gra J)^{\perp}
=\{(A^{\intercal}w, B^{\intercal}w)\mid\ w\in\vm(BA^{\intercal})+\vn(BA^{\intercal})\},$$
from which $\gra J^{*}=\gra\whg $. Then  by (i),
\begin{align*}
\gra \widetilde{G}=\gra(\whg)^*=\gra J^{**}=\gra J.
\end{align*}
\end{proof}

We are ready to apply Brezis-Browder Theorem, namely Fact~\ref{Sv:7},
 to improve Crouzeix-Anaya's characterizations
of monotonicity and maximal monotonicity of $G$ and provide a different proof.
\begin{theorem}\label{crouzeix1}Let $\whg,\wtg$ be defined in Lemma~\ref{wtgtog}.
 The following are equivalent:
\begin{enumerate}
\item
$G$ is monotone;
\item $\wtg$ is monotone;
\item $\wtg$ is maximally monotone;
\item
$\whg$ is maximally monotone;
\item  $ \dim \vp(BA^{\intercal})=p-n$, equivalently, $AB^{\intercal}+BA^{\intercal}$ has
exactly $p-n$ positive eigenvalues (counting multiplicity).
\end{enumerate}
\end{theorem}
\begin{proof}
(i)$\Leftrightarrow$(ii): Lemma~\ref{wtgtog}(iii) and Corollary~\ref{t:split}(i).

(ii)$\Leftrightarrow$(iii)$\Leftrightarrow$(iv):
Note that $\wtg=\big(\whg\big)^*$ and $\whg$ is always a monotone linear relation by Corollary~\ref{t:split}(iv). It suffices to
combine Lemma~\ref{wtgtog} and Fact~\ref{Sv:7}.

``(i)$\Rightarrow$(v)": Assume that $G$ is monotone. Then $\wtg$ is monotone by Lemma~\ref{wtgtog}(iii) and Corollary~\ref{t:split}(i).
By Lemma~\ref{wtgtog}(ii), Corollary~\ref{t:split}(iv) and Fact~\ref{Sv:7},
$\whg$ is maximally monotone, so that $\dim (\gra \whg)=p-k=n$ by Fact~\ref{TheSIAM} and Lemma~\ref{t:dimensionnochange}, thus $k=p-n$.
Note that for each eigenvalue of a symmetric matrix, its geometric multiplicity is the same as its algebraic
multiplicity \cite[page 512]{Meyer}.

``(v)$\Rightarrow$(i)": Assume that $k=p-n$. Then $\dim (\gra \whg)=p-k=n$ by Lemma~\ref{t:dimensionnochange}, so that $\whg$ is maximally monotone
by  Fact~\ref{TheSIAM}(i)(ii).
By Lemma~\ref{wtgtog}(ii) and Fact~\ref{Sv:7}, $\wtg$ is monotone, which implies that $G$ is monotone.
\end{proof}

\begin{corollary}\label{expliciteq}
Assume that $G$ is monotone. Then
\begin{align*}\gra\widetilde{G}&=\gra G+\left\{\left(\begin{matrix}B^{\intercal}\\
A^{\intercal}\end{matrix}\right)u\mid u\in \vp(BA^{\intercal}\rangle\right\}\\
&=\{(x,x^*)\mid V_g Ax+V_g Bx^*=0\}
\end{align*}
is a maximally monotone extension of $G$,
where
$$V_g=\left(\begin{matrix}
v_{p-n+1}^{\intercal}\\
v_{p-n+2}^{\intercal}\\
\vdots\\
v_{p}^{\intercal}
\end{matrix}
\right).
$$
\end{corollary}
\begin{proof}
Combine Theorem~\ref{crouzeix1} and Lemma~\ref{wtgtog}(iii) directly.
\end{proof}

A remark is in order to compare our extension with the one by Crouzeix and Anaya.
\begin{remark}
(i).
Crouzeix-Anaya \cite{CrouAna} defines
the union of monotone extension of $G$ as
$$S=\gra G+\left\{\left(\begin{matrix}B^{\intercal}\\
A^{\intercal}\end{matrix}\right)u\mid u\in K\right\},
 \text{ where $K=\{u\in \RR^n\mid\langle u, (AB^{\intercal}+BA^{\intercal})u\rangle\geq 0 \}$}.$$
Although this is the set monotonically related to $G$, it is not monotone in general
as long as
$(AB^{\intercal}+BA^{\intercal})$ has both positive eigenvalues and negative eigenvalues.
Indeed, let $(\alpha_{1}, u_{1})$ and $(\alpha_{2},u_{2})$ be eigen-pairs of $(AB^{\intercal}+BA^{\intercal})$
with $\alpha_{1}>0$ and $\alpha_{2}<0$.
We have
$$\langle u_{1}, (AB^{\intercal}+BA^{\intercal})u_{1}\rangle=\alpha_{1}\|u_{1}\|^2>0,
\quad \langle u_{2}, (AB^{\intercal}+BA^{\intercal})u_{2}\rangle=\alpha_{2}\|u_{2}\|^2<0.$$
Choose $\epsilon>0$ sufficiently small so that
$$\langle u_{1}+\epsilon u_{2}, (AB^{\intercal}+BA^{\intercal})(u_{1}+\epsilon u_{2})\rangle>0.$$
Then
$$\left(\begin{matrix}B^{\intercal}\\
A^{\intercal}\end{matrix}\right)u_{1}, \left(\begin{matrix}B^{\intercal}\\
A^{\intercal}\end{matrix}\right)(u_{1}+\epsilon u_{2})\in S.$$
However,
$$\left(\begin{matrix}B^{\intercal}\\
A^{\intercal}\end{matrix}\right)(u_{1}+\epsilon u_{2})-\left(\begin{matrix}B^{\intercal}\\
A^{\intercal}\end{matrix}\right)u_{1}=\epsilon\left(\begin{matrix}B^{\intercal}\\
A^{\intercal}\end{matrix}\right)u_{2}$$
has
$$\langle \epsilon B^{\intercal} u_{2}, \epsilon A^{\intercal}u_{2}\rangle=\epsilon^{2}\langle u_{2}, BA^{\intercal}u_{2}\rangle
=\epsilon^{2}\frac{\langle u_{2}, (AB^{\intercal}+BA^{\intercal})u_{2}\rangle}{2}<0.$$
Therefore $S$ is not monotone. By using $\vp(BA^{\intercal})\subseteq K$, we have obtained a maximally monotone extension
of $G$ .

(ii). Crouzeix and Anaya \cite{CrouAna} find the maximal monotone linear subspace extension of $G$ algorithmically by using
$u\in \widetilde{G_{k}}\setminus G_{k}$. Computationally, it is not completely clear to us how to find such an $u$.
\end{remark}

The following result extends the characterization of maximally monotone
linear relations given by Crouzeix-Anaya \cite{CrouAna}.

\begin{theorem}\label{maxchar} Let $\whg,\wtg$ be defined in Lemma~\ref{wtgtog}. The following are equivalent:
\begin{enumerate}
\item
$G$ is maximally monotone;
\item $p=n$ and $G$ is monotone;
\item $p=n$ and $AB^{\intercal}+BA^{\intercal}$ is negative semidefinite.
\item $p=n$ and $\whg$ is maximally monotone.
\end{enumerate}
\end{theorem}
\begin{proof}
(i)$\Rightarrow$(ii): Apply Lemma~\ref{MonL:1}(iii).

(ii)$\Rightarrow$(iii): Apply directly Theorem~\ref{crouzeix1}(i)(v).

(iii)$\Rightarrow$(i): Assume that $p=n$ and $(AB^{\intercal}+BA^{\intercal})$ is negative semidefinite.
Then $k=0$ and $\wtg=G$.
It follows that $\dim(\gra \whg)=p-k=n$ by Lemma~\ref{t:dimensionnochange}, so that $\whg$ is maximally monotone by Corollary~\ref{t:split}(iv) and Fact~\ref{TheSIAM}(i)(ii).
  Since $\big(\whg\big)^*=\wtg$ by Lemma~\ref{wtgtog}(ii),
Fact~\ref{Sv:7}
gives that $\wtg=G$ is maximally monotone.

(iii)$\Rightarrow$(iv): Assume that $p=n$ and $(AB^{\intercal}+BA^{\intercal})$ is negative semidefinite.
We have $k=0$ and $\dim (\gra\whg)=p-k=n-0=n$. Hence (iv) holds by Corollary~\ref{t:split}(iv) and Fact~\ref{TheSIAM}(i)(ii).

(iv)$\Rightarrow$(iii): Assume that $\whg$ is maximally monotone and $p=n$. We have
$\dim(\gra\whg)=p-k=n-k=n$ so that $k=0$. Hence $(AB^{\intercal}+BA^{\intercal})$ is negative semidefinite.
\end{proof}

We end this section with a characterization of all the maximally monotone linear subspace extensions of $G$.
\begin{theorem}\label{Genallt:1}
Let
$G$ be monotone. Then  $\widetilde{G}$ is a maximally monotone extension of $G$ if and only if there exists $N \in\RR^{p\times p}$ with rank of $n$
such that $N^{\intercal}\Id_{\lambda}N$ is negative semidefinite and
\begin{align}\gra \widetilde{G}=\{(x,x^*)\mid N^{\intercal}V^{\intercal}Ax+N^{\intercal}V^{\intercal}Bx^*=0\}.\label{Inverse:a11}\end{align}
\end{theorem}

\begin{proof}
``$\Rightarrow$'':
By Lemma~\ref{MonL:1}(i), we have
\begin{align}
\gra G^*=\{(B^{\intercal}u,-A^{\intercal}u)\mid u\in\RR^p\}.\label{Inverse:2}\end{align}
Since $\gra G\subseteq\gra \widetilde{G}$ and thus $\gra (\widetilde{G})^*$ is a subspace of
$\gra G^*$.

Thus by \eqref{Inverse:2}, there exists  a subspace $F$ of $\RR^p$ such that
\begin{align}
\gra (\widetilde{G})^*=\{(B^{\intercal}u,-A^{\intercal}u)\mid u\in F\}.\label{Inversec:3}\end{align}
By Fact~\ref{Sv:7}, Fact~\ref{TheSIAM} and Lemma~\ref{t:dimensionnochange}, we have
\begin{align}
\dim F=n.\end{align}
Thus, there exists $N\in\RR^{p\times p}$ with rank $n$ such that $\ran VN=F$ and
\begin{align}
\gra (\widetilde{G})^*=\{(B^{\intercal}VNy,-A^{\intercal}VNy)\mid  y\in \RR^p\}.\label{Inversec:4}\end{align}

As $\wtg$ is maximal monotone, $(\wtg)^*$ is maximal monotone by Fact~\ref{Sv:7}, so
$$N^{\intercal} V^{\intercal} (BA^{\intercal}+AB^{\intercal})VN \text{ is negative semidefinite.}$$
Using Fact~\ref{FaAE:1},
we have
\begin{align}N^{\intercal}\Id_{\lambda}N
=N^{\intercal}V^{\intercal}V\Id_{\lambda}N
=N^{\intercal}V^{\intercal}(AB^{\intercal}+BA^{\intercal})VN\label{Inversede:2}
\end{align}
which is negative semidefinite.
\eqref{Inverse:a11} follows from \eqref{Inversec:4} by Lemma~\ref{FaAEL:2}.

``$\Leftarrow$'':
By Lemma~\ref{FaAEL:2}, we have
\begin{align}
\gra (\widetilde{G})^*=\{(B^{\intercal}VNu,-A^{\intercal}VNu)\mid u\in\RR^p\}.\label{Inversed:2}\end{align}
Observe that $(\widetilde{G})^*$ is monotone because
$N^{\intercal}V^{\intercal}(AB^{\intercal}+BA^{\intercal})VN=N^{\intercal}\Id_{\lambda}N$ is negative semidefinite
by Fact~\ref{FaAE:1} and the assumption.  As $\rank (VN)=n$, it follows from
\eqref{Inversed:2} and Lemma~\ref{t:dimensionnochange} that $\dim\gra (\widetilde{G})^*=n$.
Therefore  $(\widetilde{G})^*$ is maximally monotone by Fact~\ref{TheSIAM}. Applying Fact~\ref{Sv:7} for $ T=(\widetilde{G})^*$
 yields that $\widetilde{G}=(\widetilde{G})^{**}$ is maximally monotone.
\end{proof}

From the above proof, we see that to find a maximal monotone extension
extension of $G$ one essentially need to find subspace $F\subseteq \RR^p$ such that $\dim F=n$ and
$AB^{\intercal}+BA^{\intercal}$ is negative semidefinite on $F$.
If $F=\ran M$ and $M\in\RR^{p\times p}$  with $\rank M=n$, one can let
$N=V^{\intercal}M$. The maximal monotone linear subspace extension of $G$ is
$$\wtg=\{(x,x^*)\mid M^{\intercal}Ax+M^{\intercal}Bx^*=0\}.$$
In Corollary~\ref{expliciteq}, one can choose $M=\big(\underbrace{0\ 0\ \cdots 0}_{n}\ \ \ v_{p-n+1}\ \cdots \ v_{p}\big)$.

\begin{corollary}\label{t:Mform}
Let
$G$ be monotone. Then  $\widetilde{G}$ is a maximally monotone extension of $G$ if and only if there exists $M\in\RR^{p\times p}$ with rank of $n$
such that $M^{\intercal}(AB^{\intercal}+BA^{\intercal})M$ is negative semidefinite and
\begin{align}\gra \widetilde{G}=\{(x,x^*)\mid M^{\intercal} Ax+M^{\intercal}Bx^*=0\}.\label{Inverse:a3}\end{align}
\end{corollary}

Note that $G$ may have different representations in terms of $A,B$. The maximal monotone extension of $\wtg$ given in Theorem~\ref{Genallt:1}
and Corollary~\ref{expliciteq}
relies on $A,B$ matrices
and $N$.
This might leads different maximal monotone extensions, see Section~\ref{Exam:sec}.


\section{Minty parameterizations}\label{Minty}
Although $G$ is set-valued in general, when $G$ is monotone it has a beautiful Minty parametrization in terms of $A,B$, which is what we are going to show
in this section.

\begin{lemma}
The linear relation $G$ is monotone if and only if
\begin{align}
&\|y\|^2-\|y^*\|^2\geq 0, \text{ whenever }\label{objective}\\
& (A+B)y+(B-A)y^*=0.\label{constraint}
\end{align}
Consequently, if $G$ is monotone then the $p\times n$ matrix $B-A$ must have full column rank, namely $n$.
\end{lemma}
\begin{proof}
Define
$$P=\left(\begin{matrix}
0 & \Id\\
\Id & 0
\end{matrix}
\right).
$$
It is easy to see that $G$ is monotone if and only if
$$ \langle (x,x^*), P\left(\begin{matrix} x\\
x^*\end{matrix}
\right)\rangle\geq 0,$$
whenever
$Ax+Bx^*=0$.
Define the orthogonal matrix
$$Q=\frac{1}{\sqrt{2}}\left(\begin{matrix}
\Id & -\Id\\
\Id & \Id
\end{matrix}
\right)
$$
and put
$$\left(\begin{matrix} x\\
x^*
\end{matrix}
\right)=Q\left(\begin{matrix} y\\
y^*
\end{matrix}
\right).$$
Then $G$ is monotone if and only if
\begin{align}
&\|y\|^2-\|y^*\|^2\geq 0, \text{ whenever }\label{objectiveold}\\
& (A+B)y+(B-A)y^*=0. \label{constraintold}
\end{align}
If $(B-A)$ does not have full column rank, then there exists $y^*\neq 0$ such that $(B-A)y^*=0$.
Then $(0,y^*)$ satisfies \eqref{constraintold} but \eqref{objectiveold} fails.
Therefore, $B-A$ has to be full column rank.
\end{proof}
\begin{theorem}[Minty parametrization] \label{t:param}
Assume that $G$ is a monotone operator. Then  $(x,x^*)\in\gra G$
if and only if
\begin{align}
x &=\frac{1}{2}[\Id +(B-A)^{\dagger}(B+A)]y\label{xxall1}\\
x^* &=\frac{1}{2}[\Id-(B-A)^{\dagger}(B+A)]y\label{xxall2}
\end{align}
for $y=x+x^*\in\ran(\Id+G)$. Here the Moore-Penrose inverse
$(B-A)^{\dagger}=[(B-A)^{\intercal}(B-A)]^{-1}(B-A)^{\intercal}$.
In particular, when $G$ is maximally monotone, we have
$$\gra G=\{((B-A)^{-1}By, -(B-A)^{-1}Ay)\mid y\in \RR^n\}.$$
\end{theorem}
\begin{proof}
As $(B-A)$ is full column rank, $(B-A)^{\intercal}(B-A)$ is invertible. It follows from \eqref{constraint} that
$(B-A)^{\intercal}(A+B)y+(B-A)^{\intercal}(B-A)y^*=0$ so that
$$y^*=-((B-A)^{\intercal}(B-A))^{-1}(B-A)^{\intercal}(A+B)y=-(B-A)^{\dagger}(A+B)y.$$
Then
\begin{align*}
x&=\frac{1}{\sqrt{2}}(y-y^*)=\frac{1}{\sqrt{2}}[\Id +(B-A)^{\dagger}(B+A)]y\\
x^*&=\frac{1}{\sqrt{2}}(y+y^*)=\frac{1}{\sqrt{2}}[\Id-(B-A)^{\dagger}(B+A)]y
\end{align*}
where $y=\frac{x+x^*}{\sqrt{2}}$ with $(x,x^*)\in \gra G$. Since $\ran (\Id+G)$ is a subspace, we
have
\begin{align*}
x&=\frac{1}{2}[\Id +(B-A)^{\dagger}(B+A)]\tilde{y}\\
x^*&=\frac{1}{2}[\Id-(B-A)^{\dagger}(B+A)]\tilde{y}
\end{align*}
with $\tilde{y}=x+x^*\in\ran(\Id+G)$.

If $G$ is maximally monotone, then $p=n$ by Theorem~\ref{maxchar} and hence $B-A$ is invertible, thus $(B-A)^{\dagger}=(B-A)^{-1}$.
Moreover, $\ran (G+\Id)=\RR^n$.
Then \eqref{xxall1} and \eqref{xxall2} transpire to
\begin{align}
x &=\frac{1}{2}(B-A)^{-1}[B-A+(B+A)]y=(B-A)^{-1}By\\
x^* &=\frac{1}{2}(B-A)^{-1}[(B-A)-(B+A)]y=-(B-A)^{-1}Ay
\end{align}
for $y\in\RR^n$.
\end{proof}

\begin{remark} See Lemma~\ref{AuSe:1} for $\ran (G+\Id)$.
Note that as $G$ is a monotone linear relation,
the mapping
$$z\mapsto ((G+\Id)^{-1}, \Id-(G+\Id)^{-1})(z)$$
is bijective and linear from $\ran (G+\Id)$ to $\gra G$,
therefore
$\dim(\ran (G+\Id))=\dim (\gra G)$.
\end{remark}

\begin{corollary}
Let $G$ be a monotone operator. Then $\wtg$ define in Corollary~\ref{expliciteq}, the maximally monotone extension of $G$, has its
Minty parametrization given by
$$\gra \wtg=\{((V_g B-V_g A)^{-1}V_g By, -(V_g B-V_g A)^{-1}V_g Ay)\mid y\in \RR^n\}$$
where $V_g$ is given as in Corollary~\ref{expliciteq}.
\end{corollary}
\begin{proof}
Since $\rank (V_g)=n$ and $\rank(A\ B)=p$, by Lemma~\ref{t:dimensionnochange}\eqref{t:dim+}, $\rank (V_g A\ V_g B)=n$. Apply
 Corollary~\ref{expliciteq} and Theorem~\ref{t:param} directly.
\end{proof}

\begin{corollary}
When $G$ is maximally monotone,
$$\dom G =(B-A)^{-1}(\ran B), \quad \ran G= (B-A)^{-1}(\ran A).$$
\end{corollary}
Recall that $T:\RR^n\rightarrow\RR^n$ is \emph{firmly nonexpansive} if
$$\|Tx-Ty\|^2\leq\scal{Tx-Ty}{x-y}\quad \forall \ x,y\in\dom T.$$
In terms of matrices
\begin{corollary}
Suppose that $p=n$, $AB^{\intercal}+BA^{\intercal}$ is negative semidefinite.  Then
$(B-A)^{-1}B$ and $-(B-A)^{-1}A$ are firmly nonexpansive.
\end{corollary}
\begin{proof} By Theorem~\ref{maxchar}, $G$ is maximal monotone. Theorem~\ref{t:param} gives that
$$(B-A)^{-1}B=(\Id+G)^{-1}, \quad -(B-A)^{-1}A=(\Id+G^{-1})^{-1}.$$
Being resolvent of monotone operators $G, G^{-1}$, they are firmly
nonexpansive, see \cite{BC2011,EcBe} or \cite[Fact~2.5]{bw10}.
\end{proof}

\section{Maximally monotone extensions with the same domain or the same range}\label{DRan:1}

The purpose of this section is to find maximal monotone linear subspace extensions of $G$
which keep either $\dom G$ or $\ran G$ unchanged.
For a closed convex set $S\subseteq \RR^n$, let $N_{S}$ denote its normal cone mapping.
\begin{proposition}\label{domsameextension}
Assume that $T:\RR^n\To\RR^n$ is a monotone linear relation. Then
\begin{enumerate}
\item
$T_{1}=T+N_{\dom T}$, i.e.,
$$x\mapsto T_{1}x=\begin{cases}
 Tx+(\dom T)^{\perp} &\text{ if $x\in \dom T$} \\
 \emptyset & \text{ otherwise}
 \end{cases}
 $$ is
 maximally monotone. In particular, $\dom T_{1}=\dom T$.
 \item
 $T_{2}=(T^{-1}+N_{\ran T})^{-1}$ is a maximally monotone extension
 of $T$ and $\ran T_{2}=\ran T$.
 \end{enumerate}
\end{proposition}

\begin{proof}
(i):\ Since $0\in T0\subseteq (\dom T)^{\perp}$ by\cite[Proposition~2.2(i)]{bwy09}, we have
$T_{1}0=T0+(\dom T)^{\perp}=(\dom T)^{\perp}$ so that
$\dom T_{1}=\dom T= (T_{1}0)^{\perp}.$
Hence $T_{1}$ is maximally monotone by Fact~\ref{TheSIAM}.

(ii):\ Apply (i) to $T^{-1}$ to see that
$T^{-1}+N_{\ran T}$ is a maximally monotone extension of $T^{-1}$ with $\dom (T^{-1}+N_{\ran T})=\ran T$.
 Therefore,
$T_{2}$ is a maximally monotone extension of $T$ with $\ran T_{2}=\ran T$.
\end{proof}

Since
$$\gra G=\{(x,x^*)\mid\ Ax+Bx^*=0\}$$
we can use Gaussian elimination to reduce (A\ B) to row echelon form. Then back substitution
to solve basic variables in terms of the free variables, see \cite[page 61]{Meyer}. Row-echelon form gives
$$\left(\begin{matrix}
x\\
x^*
\end{matrix}
\right)=h_{1}y_{1}+\cdots +h_{2n-p}y_{2n-p}=\left(\begin{matrix}
C\\
D
\end{matrix}
\right)y
$$
where $y\in \RR^{2n-p}$ and
$$\left(\begin{matrix}
C\\
D
\end{matrix}
\right)=(h_{1},\ldots, h_{2n-p})$$
with $C,D$ being $n\times (2n-p)$ matrices. Therefore,
\begin{equation}\label{t:range}
\gra G=\left\{\left(\begin{matrix}
Cy\\
Dy
\end{matrix}
\right)\mid\ y\in \RR^{2n-p}\right\}.
\end{equation}
Define
\begin{equation}\label{t:rangeform}
\gra E_{1}=\left\{\left(\begin{matrix}
Cy\\
Dy
\end{matrix}
\right)+\left(\begin{matrix}
0\\
(\ran C)^{\perp}\end{matrix}
\right)\mid\ y\in \RR^{2n-p}\right\}.
\end{equation}

\begin{equation}\label{t:rangeform2}
\gra E_{2}=\left\{\left(\begin{matrix}
Cy\\
Dy
\end{matrix}
\right)+\left(\begin{matrix}
(\ran D)^{\perp}\\
0\end{matrix}
\right)\mid\ y\in \RR^{2n-p}\right\}.
\end{equation}

\begin{theorem}\label{t:keepdr}
\begin{enumerate}
\item
$E_{1}$ is a maximally monotone extension of $G$ with
$\dom E_{1}=\dom G$. Moreover,
\begin{equation}\label{rangeform}
\gra E_{1}=\ran \left(\begin{matrix}
C\\
D
\end{matrix}
\right)+\left(\begin{matrix}
0\\
(\ran C)^{\perp}\end{matrix}
\right)=\ran \left(\begin{matrix}
C\\
D
\end{matrix}
\right)+\left(\begin{matrix}
0\\
\ker C^{\intercal}\end{matrix}
\right).
\end{equation}
\item
$E_{2}$ is a maximally monotone extension of $G$ with
$\ran E_{2}=\ran G$. Moreover,
\begin{equation}\label{domainform}
\gra E_{2}=\ran \left(\begin{matrix}
C\\
D
\end{matrix}
\right)+\left(\begin{matrix}
(\ran D)^{\perp}\\
0\end{matrix}
\right)=\ran \left(\begin{matrix}
C\\
D
\end{matrix}
\right)+\left(\begin{matrix}
\ker D^{\intercal}\\
0\end{matrix}
\right).
\end{equation}
\end{enumerate}
\end{theorem}
\begin{proof}
(i):\
Note that $\dom G=\ran C$. The maximal monotonicity follows from Proposition~\ref{domsameextension}.
\eqref{rangeform} follows from \eqref{t:rangeform} and that $(\ran C)^{\perp}=\ker C^{\intercal}$
\cite[page 405]{Meyer}.

(ii): Apply (i) to $G^{-1}$, i.e.,
\begin{equation}
\gra G^{-1}=\left\{\left(\begin{matrix}
Dy\\
Cy
\end{matrix}
\right)\mid\ y\in \RR^{2n-p}\right\}
\end{equation}
and followed by taking the set-valued inverse.
\end{proof}

Apparently, both extensions $E_{1},E_{2}$ rely on $\gra G$, $\dom G, \ran G$, not on the
$A,B$. In this
sense, $E_{1},E_{2}$ are intrinsic maximal monotone linear subspace extensions.

\begin{remark}
Theorem~\ref{t:keepdr} is much easier to use than Corollary~\ref{t:Mform}
when $G$ is written in the form of \eqref{t:range}. Indeed, it is not hard to check that
\begin{align}\gra (E_1^*)=\{(B^{\intercal}u, -A^{\intercal}u)\mid B^{\intercal}u\in\dom G, u\in \RR^p\}.\label{Rem:a2}\end{align}
\begin{align}\gra (E_2^*)=\{B^{\intercal}u, -A^{\intercal}u)\mid A^{\intercal}u\in\ran G, u\in\RR^p\}.\label{Rem:a3}\end{align}
According to Fact~\ref{Sv:7}, $E_{i}^*$ is maximal monotone and $\dim E_{i}^*=n$. This implies that
$$\dim \{u\in \RR^p\mid B^{\intercal}u\in\dom G\}=n,\quad \quad \dim\{u\in\RR^{p}\mid A^{\intercal}u\in\ran G\}=n.$$
Let $M_{i}\in\RR^{p\times p}$ with $\rank M=n$ and
\begin{equation}\label{amyone}
\{u\in \RR^p\mid B^{\intercal}u\in\dom G\}=\ran M_{1},
\end{equation}
\begin{equation}\label{amytwo}
\{u\in\RR^p\mid A^{\intercal}u\in\ran G\}=\ran M_{2}.
\end{equation}
Corollary~\ref{t:Mform} shows that
$$\gra E_{i}=\{(x,x^*)\mid M_{i}^{\intercal} Ax+M_{i}^{\intercal} Bx^*=0\}.$$
However, finding $M_{i}$ from \eqref{amyone} and \eqref{amytwo} may not be easy as it seems.
\end{remark}

\section{Examples}\label{Exam:sec}
In the final section,  we illustrate our maximally monotone extensions by considering
three examples. In particular, they show that maximal monotone extensions $\wtg$ rely on
the representation of $G$ in terms of $A,B$ and choices of $N$ we shall use. However, the maximal monotone extensions
$E_{i}$ are intrinsic, only depending on $\gra G$.

\begin{example}\label{february2} Consider
$$\gra G=\left\{\left(\begin{matrix}
x\\
x^*
\end{matrix}
\right)\mid \left(\begin{matrix}
\Id\\
0
\end{matrix}\right)x+\left(\begin{matrix}
0\\
C
\end{matrix}\right)x^*=0, x,x^*\in \RR^n\right\}$$
where $C$ is a $n\times n$ symmetric, positive definite matrix. Clearly,
$$\gra G=\left\{\left(\begin{matrix}
0\\
0\end{matrix}
\right)
\right\}.$$
We have
\begin{enumerate}
\item For every $\alpha\in\left[-1,1\right], \wtg_{\alpha} $ defined by \begin{align*}
\gra\wtg_{\alpha} =\begin{cases} \left\{(0,\RR^n) \right\},\ &\text{if}\ \alpha=1;\\
 \left\{(x,\tfrac{1+\alpha}{1-\alpha}C^{-1}x)\mid x\in\RR^n \right\},\ &\text{otherwise}\end{cases}
\end{align*}
 is a maximally monotone linear extension of $G$.

\item $ E_{1}=\wtg_1$ and
$E_{2}=\wtg_{-1}.$
\end{enumerate}
\end{example}
\begin{proof}
(i): To find $\wtg_{\alpha} $, we need eigenvectors of
$$\bA=\left(\begin{matrix}
\Id\\
0
\end{matrix}\right)
(0 \ C^{\intercal})+\left(\begin{matrix}
0\\
C
\end{matrix}
\right) (\Id\ 0)
=\left(\begin{matrix}
0 & C \\
C & 0
\end{matrix}
\right)$$
Counting multiplicity, the positive definite matrix $C$ has eigen-pairs $(\lambda_{i},w_{i})$
($i=1,\ldots, n$)
such that $\lambda_{i}>0, \|w_{i}\|=1$ and $\langle w_{i},w_{j}\rangle =0$ for
$i\neq j$.
As such, the matrix
$\bA$
has $2n$
eigen-pairs, namely
$$
(\lambda_{i}, \left(\begin{matrix}
w_{i}\\
w_{i}
\end{matrix}
\right))
$$
and
$$(-\lambda_{i}, \left(\begin{matrix}
w_{i}\\
-w_{i}
\end{matrix}
\right))
$$
with $i=1,\ldots, n$. Put $W=[w_{1}\ \cdots \ w_{n}]$ and write
$$V=\begin{pmatrix}
W & W\\
W & -W
\end{pmatrix}.
$$
Then
$$W^{\intercal}CW=
D=\text{diag}(\lambda_{1},\lambda_{2},\cdots, \lambda_{n})
$$
In Theorem~\ref{Genallt:1}, take
 \begin{align*}
 N_{\alpha}=
 \begin{pmatrix}
 \mathbf{0}& \alpha\Id\\
 \mathbf{0} & \Id \\
\end{pmatrix}.
 \end{align*}
 We have  $\rank N_{\alpha}=n$,
 \begin{align*}
 N_{\alpha}^{\intercal} \Id_{\lambda} N_{\alpha}=
 \begin{pmatrix}
 \mathbf{0} & \mathbf{0}\\
 \mathbf{0} & (\alpha^2-1)W^{\intercal}CW
\end{pmatrix}= \begin{pmatrix}
 \mathbf{0} & \mathbf{0}\\
 \mathbf{0} & (\alpha^2-1)D
\end{pmatrix}
 \end{align*}
 being  negative semidefinite,
and
$$VN_{\alpha}= \begin{pmatrix}
0 & (1+\alpha)W\\
0 & (\alpha -1) W
\end{pmatrix}.
$$
Then by Theorem~\ref{Genallt:1}, we have an maximally monotone linear extension $\wtg_{\alpha}$ given by
\begin{align*}
\gra \wtg_{\alpha} &=\left\{(x,x^*)\mid \left(\begin{matrix}0\\ (1+\alpha)W^{\intercal}x+(\alpha-1)W^{\intercal}Cx^*
 \end{matrix}\right)=0\right\}\\
 &=\left\{(x,x^*)\mid  (1+\alpha)x+(\alpha-1)Cx^*=0\right\}\\
 &=\begin{cases} \left\{(0,\RR^n) \right\},\ &\text{if}\ \alpha=1;\\
 \left\{(x,\tfrac{1+\alpha}{1-\alpha}C^{-1}x)\mid x\in\RR^n \right\},\ &\text{otherwise}\end{cases},
\end{align*}
Hence we get the result as desire.

(ii): It is  immediate from Theorem~\ref{t:keepdr} and (i).
\end{proof}

\begin{example}\label{february1}
Consider
$$\gra G=\left\{\left(\begin{matrix}
x\\
x^*
\end{matrix}
\right)\mid \left(\begin{matrix}
-1 & 0\\
0 & 0 \\
0 & -1
\end{matrix}\right)\left(\begin{matrix}
x_{1}\\
x_{2}
\end{matrix}
\right)+\left(\begin{matrix}
1 & 0\\
0 & 1\\
0 & 1
\end{matrix}\right)\left(\begin{matrix}
x_{1}^*\\
x_{2}^*
\end{matrix}
\right)=0, x_i,x_{i}^*\in \RR\right\}.$$
Then
\begin{enumerate}
\item
\begin{align*}\wtg_1=\left(\begin{matrix}
1 & 0\\
0 &\frac{-1+\sqrt{2}}{2-\sqrt{2}}
\end{matrix}\right),\quad \wtg_2=\begin{pmatrix}
1&\frac{2}{5}\\
0&\frac{\sqrt{2}}{10}\end{pmatrix}
\end{align*}
are the maximally monotone extensions of $G$.
\item $$E_{1}(x_{1},0)=(x_{1},\RR)\quad \forall x_{1}\in \RR.$$

\item $$
E_{2}(x_{1},y)=(x_{1},0)\quad \forall x_{1}, y\in\RR.$$
\end{enumerate}
\end{example}

\begin{proof}
We have
$$\gra G=\left\{\left(\begin{matrix}
x_{1}\\
0\\
x_{1}\\
0
\end{matrix}
\right)\mid x_{1}\in \RR\right\}$$
is monotone.
Since $\dim G=1$, $G$ is not maximally monotone by Fact~\ref{TheSIAM}.

The matrix
$$AB^{\intercal}+BA^{\intercal}=\left(\begin{matrix}
-2 & 0 & 0\\
0 & 0 & -1\\
0 & -1 & -2
\end{matrix}
\right)
$$
has positive eigenvalue $-1+\sqrt{2}$ with eigenvector
$$u=\left(\begin{matrix}
0\\
1\\
1-\sqrt{2}
\end{matrix}
\right)\quad \text{ so that } \left(\begin{matrix}
B^{\intercal}\\
A^{\intercal}
\end{matrix}
\right)u=\left(\begin{matrix}
0\\
2-\sqrt{2}\\
0\\
-1+\sqrt{2}
\end{matrix}
\right).$$
Then by Corollary~\ref{expliciteq},
$\gra \wtg_1=
$
$$\left\{\left(\begin{matrix}
x_{1}\\
0\\
x_{1}\\
0
\end{matrix}
\right)\mid x_{1}\in \RR\right\}+\left\{\left(\begin{matrix}
0\\
2-\sqrt{2}\\
0\\
-1+\sqrt{2}
\end{matrix}
\right)x_{2}\mid x_{2}\in \RR\right\}=\left\{\left(\begin{matrix}
x_{1}\\
(2-\sqrt{2})x_{2}\\
x_{1}\\
(-1+\sqrt{2})x_{2}
\end{matrix}
\right)\mid x_{1}, x_{2}\in \RR\right\}.$$
Therefore,
$$\wtg_1=\left(\begin{matrix}
1 & 0\\
0 &\frac{-1+\sqrt{2}}{2-\sqrt{2}}
\end{matrix}
\right).
$$
We have
\begin{align}\label{back1}
\Id_{\lambda}=\begin{pmatrix}
-1+\sqrt{2}&0&0\\
0&-1-\sqrt{2}&0 \\
0&0&-2\\
\end{pmatrix},\quad V=\begin{pmatrix}
0&0&1\\
-\frac{1}{-1+\sqrt{2}}&-\frac{1}{-1-\sqrt{2}}&0 \\
1&1&0\\
\end{pmatrix}.
\end{align}

Take
\begin{align}\label{back2}
N=\begin{pmatrix}
0&-1 & 1\\
0 &2&-1\\
0&1&1
\end{pmatrix}.
\end{align}
We have $\rank N=2$ and
\begin{align}
N^{\intercal}\Id_{\lambda}N=\left(\begin{matrix}
0&0 & 0\\
0 &-7-3\sqrt{2}&1+\sqrt{2}\\
0&1+\sqrt{2}&-4
\end{matrix}
\right),
\end{align}
being negative semidefinite.

By Theorem~\ref{Genallt:1}, with $V, N$ given in \eqref{back1} and \eqref{back2},
we use the NullSpace command in Maple to solve
$$(VN)^{\intercal}Ax+(VN)^{\intercal}Bx^*=0,$$
and get
\begin{align*}
\gra \wtg_2=\spand\left\{\ \begin{pmatrix}
1\\
0\\
1\\
0\\
\end{pmatrix}\ \begin{pmatrix}
-2\sqrt{2}\\
5\sqrt{2} \\
0\\
1\\
\end{pmatrix}\right\}.
\end{align*}
Thus
$\wtg_2=\begin{pmatrix}
1&-2\sqrt{2}\\
0&5\sqrt{2}
\end{pmatrix}^{-1}=\begin{pmatrix}
1&\frac{2}{5}\\
0&\frac{\sqrt{2}}{10}
\end{pmatrix}$.

On the other hand,
$$\gra E_{1}=\left\{\left(\begin{matrix}
x_{1}\\
0\\
x_{1}\\
0
\end{matrix}
\right)\mid x_{1}\in \RR\right\}+\left(\begin{matrix}
0\\
0\\
0\\
\RR
\end{matrix}
\right)
=\left\{\left(\begin{matrix}
x_{1}\\
0\\
x_{1}\\
\RR
\end{matrix}
\right)\mid x_{1}\in \RR\right\}
$$
 gives
$$E_{1}(x_{1},0)=(x_{1},\RR)\quad \forall x_{1}\in \RR.$$
And $$\gra E_{2}
=\left\{\left(\begin{matrix}
x_{1}\\
\RR\\
x_{1}\\
0
\end{matrix}
\right)\mid x_{1}\in \RR\right\}.
$$
gives $$
E_{2}(x_{1},y)=(x_{1},0)\quad \forall x_{1}, y\in\RR.$$
\end{proof}

In \cite{BW}, the authors use autoconjugates to find maximally monotone extensions
of monotone operators. In general, it is not clear whether the maximally monotone extensions
of a linear relation is still a linear relation. As both monotone operators in Examples~\ref{february1} and Examples~\ref{february2}
are subset of $\{(x,x)\mid\ x\in \RR^n\}$, \cite[Example~5.10]{BW} shows that the maximally monotone extension
obtained by autoconjugate must be $\Id$, which are different from the ones given here.

\begin{example} Set  $\gra G=\{(x,x^*)\mid Ax+Bx^*=0\}$ where
\begin{align}\label{t:abneeded}
A=\begin{pmatrix}
1&1 \\
2&0\\
3&1\\
\end{pmatrix},
\, B=\begin{pmatrix}
1&5 \\
1&7\\
0&2
\end{pmatrix},\text{thus}\, (A\ B)=\begin{pmatrix}
1&1&1&5 \\
2&0&1&7\\
3&1&0&2
\end{pmatrix}.
\end{align}
Then
\begin{align*}
\widetilde{G}_{1}=\begin{pmatrix}\frac{-117+17\sqrt{201}}{2(-1+\sqrt{201})}
&\frac{-107+7\sqrt{201}}{2(-1+\sqrt{201})}\\
-\frac{-23+3\sqrt{201}}{2(-1+\sqrt{201})}&-\frac{-21+\sqrt{201}}{2(-1+\sqrt{201})}\\
\end{pmatrix},\quad\widetilde{G}_{2}=\begin{pmatrix}\frac{33}{4}-\frac{\sqrt{201}}{6}
&\frac{13}{4}-\frac{\sqrt{201}}{6}\\
-\frac{29}{20}+\frac{\sqrt{201}}{30}&-\frac{9}{20}+\frac{\sqrt{201}}{30}\\
\end{pmatrix}.\end{align*}
are two maximally monotone
linear extensions of $G$.

Moreover,
\begin{align*}
\gra E_{1}&=\left\{\begin{pmatrix}
-1\\
1\\
-5\\
1
\end{pmatrix}x_1 + \begin{pmatrix}
0\\
0\\
1\\
1
\end{pmatrix}x_2 \mid x_1, x_2\in\RR\right\},\quad
\gra E_{2}=\left\{\begin{pmatrix}
-1\\
1\\
-5\\
1
\end{pmatrix}x_1 +\begin{pmatrix}
1\\
5\\
0\\
0
\end{pmatrix} x_2\mid x_1,x_2\in\RR\right\}.
\end{align*}

\begin{proof}
We have $\rank (A\ B)=3 $ and
\begin{align}\label{t:eigenvs}
\Id_{\lambda}=\begin{pmatrix}
13+\sqrt{201}&0&0\\
0&-6&0 \\
0&0&13-\sqrt{201}\\
\end{pmatrix},\quad V=\begin{pmatrix}
\tfrac{20}{1+\sqrt{201}}&0&\tfrac{20}{1-\sqrt{201}}\\
1&-1&1 \\
1&1&1\\
\end{pmatrix},
\end{align}
and
\begin{align}\label{t:vgneeded} V_g=\begin{pmatrix}
0&-1&1\\
\tfrac{20}{1-\sqrt{201}}&1&1
\end{pmatrix}.
\end{align}

Clearly, here $p=3, n=2$ and
 $AB^{\intercal}+BA^{\intercal}$ has exactly $p-n=3-2=1$
positive  eigenvalue.
By Theorem~\ref{crouzeix1}(i)(v), $G$ is monotone.

Since $AB^{\intercal}+BA^{\intercal}$ is not negative semidefinite, by Theorem~\ref{maxchar}(i)(iii),
$G$ is not maximally monotone.   With $V_{g}$ given in \eqref{t:vgneeded} and $A, B$ in \eqref{t:abneeded}, use the NullSpace command
in maple to solve
$V_{g}Ax+V_{g}Bx^*=0$
and obtain $\widetilde{G}_{1}$ defined by
\begin{align*}\gra\widetilde{G}_{1}=\spand\left\{
\begin{pmatrix}
-\frac{-21+\sqrt{201}}{2(-1+\sqrt{201})}\\
\frac{-23+3\sqrt{201}}{2(-1+\sqrt{201})} \\
1\\
0\\
\end{pmatrix}, \begin{pmatrix}
-\frac{-107+7\sqrt{201}}{2(-1+\sqrt{201})}\\
\frac{-117+17\sqrt{201}}{2(-1+\sqrt{201})} \\
0\\
1\\
\end{pmatrix}\right\}.
\end{align*}
By  Corollary~\ref{expliciteq}, $\wtg_{1}$ is a  maximally monotone
linear subspace extension of $G$.
Then
\begin{align*}\widetilde{G}_{1}=\begin{pmatrix}
-\frac{-21+\sqrt{201}}{2(-1+\sqrt{201})}&-\frac{-107+7\sqrt{201}}{2(-1+\sqrt{201})}\\
\frac{-23+3\sqrt{201}}{2(-1+\sqrt{201})}&\frac{-117+17\sqrt{201}}{2(-1+\sqrt{201})}\\
\end{pmatrix}^{-1}=\begin{pmatrix}\frac{-117+17\sqrt{201}}{2(-1+\sqrt{201})}
&\frac{-107+7\sqrt{201}}{2(-1+\sqrt{201})}\\
-\frac{-23+3\sqrt{201}}{2(-1+\sqrt{201})}&-\frac{-21+\sqrt{201}}{2(-1+\sqrt{201})}\\
\end{pmatrix}.\end{align*}

Let $N$ be defined by
\begin{align}\label{t:second}
N=\begin{pmatrix}
0&0&\frac{1}{5}\\
0&1&0 \\
0&0&1\\
\end{pmatrix}.
\end{align}
Then $\rank N=2$ and
\begin{align*}
 N^{\intercal}\Id_{\lambda}N=\begin{pmatrix}
0&0&0\\
0&-6&0 \\
0&0&\frac{338-24\sqrt{201}}{25}\\
\end{pmatrix}.
\end{align*}
 is negative semidefinite.

With $N$ in \eqref{t:second}, $A, B$ in \eqref{t:abneeded} and $V$ in \eqref{t:eigenvs},  use the NullSpace command in maple to solve
$(VN)^{\intercal}Ax+(VN)^{\intercal}Bx^*=0$. By Theorem~\ref{Genallt:1}, we get a  maximally monotone
linear extension of $G$, $\widetilde{G}_{2}$, defined by
\begin{align*}\widetilde{G}_{2}=\begin{pmatrix}
-\frac{9}{20}+\frac{\sqrt{201}}{30}&-\frac{13}{4}+\frac{\sqrt{201}}{6}\\
\frac{29}{20}-\frac{\sqrt{201}}{30}&\frac{33}{4}-\frac{\sqrt{201}}{6}\\
\end{pmatrix}^{-1}=\begin{pmatrix}\frac{33}{4}-\frac{\sqrt{201}}{6}
&\frac{13}{4}-\frac{\sqrt{201}}{6}\\
-\frac{29}{20}+\frac{\sqrt{201}}{30}&-\frac{9}{20}+\frac{\sqrt{201}}{30}\\
\end{pmatrix}.
\end{align*}

To find $E_1$ and $E_2$, using the LinearSolve command in Maple,  we  get $\gra G=\ran\begin{pmatrix}C\\ D\end{pmatrix}$, where
\begin{align*}C=\begin{pmatrix}
-1\\
1\\
\end{pmatrix},\quad D=\begin{pmatrix}
-5\\
1\\
\end{pmatrix}.
\end{align*}
It follows from Theorem~\ref{t:keepdr} that
\begin{align*}
\gra E_{1}&=\left\{\begin{pmatrix}
-1\\
1\\
-5\\
1
\end{pmatrix}x_1 + \begin{pmatrix}
0\\
0\\
1\\
1
\end{pmatrix}x_2 \mid x_1, x_2\in\RR\right\},
\end{align*}
\begin{align*}
\gra E_{2}=\left\{\begin{pmatrix}
-1\\
1\\
-5\\
1
\end{pmatrix}x_1 +\begin{pmatrix}
1\\
5\\
0\\
0
\end{pmatrix} x_2\mid x_1,x_2\in\RR\right\}.
\end{align*}
\end{proof}
\end{example}

\section*{Acknowledgments}
The authors thank Dr.\ Heinz Bauschke for bring their attentions of
\cite{CrouAna} and
many valuable discussions.
Xianfu Wang was partially supported by the Natural Sciences and
Engineering Research Council of Canada.

\end{document}